\input amstex
\documentstyle {amsppt}
\input epsf
\document
\NoBlackBoxes
\documentstyle {amsppt}
\magnification=1200
\define \cl {\Cal L_{S^1}}
\define \cg {\Cal L_{\Gamma}}
\define \Ga {\Gamma}
\define \la {\lambda}
\define \pa {\partial}
\define \bR {\Bbb R}
\define \bC {\Bbb C}
\def \bZ{\Bbb Z}
\define \tl {tree-like }
\define \tln { tree-like normalized }
\define \tlc { tree-like curve }
\def \ga {\gamma}
\def \part {\partial}
\def \pr {\prime}

\leftline{\eightit To  Vladimir Igorevich Arnold}
\leftline{\eightit with love and admiration}

\vskip 3cm

\topmatter
\title
Tree-like curves and their number of inflection points
\endtitle

\author  B.~Shapiro
\endauthor
\affil
 Department of Mathematics, University of Stockholm, S-10691, Sweden,
{\tt shapiro\@matematik.su.se}\\
\endaffil
\rightheadtext {Tree-like curves and their inflection points}
\leftheadtext {B.~Shapiro}

\abstract
In this short note we give a criterion when a planar
tree-like curve, i.e. a generic curve in $\Bbb R^2$ each double point of
which cuts it into two disjoint parts, can be send by a diffeomorphism of
$\Bbb R^2$ onto a curve with no inflection points. We also present
  some upper and lower bounds for the minimal number of inflection
points  on such curves unremovable by diffeomorphisms of $\bR^2$.
\endabstract

\endkeywords

\endtopmatter

\heading \S 1. Introduction \endheading

This paper provides a partial answer to the following question posed to the
author by V.Arnold in June 95. Given a generic immersion $c:S^1\to \bR^2$
  (i.e. with double points only)  let $\sharp_{inf}(c)$ denote the number of
inflection points on $c$ (assumed finite) and let $[c]$ denote the class of
 $c$,  i.e. the connected component in the space of generic immersions of $S^1$
 to $\bR^2$ containing $c$. Finally, let $\sharp_{inf}[c]=\min_{c^\prime\in
[c]}\sharp_{inf}(c^\prime).$
\vskip 2truemm

{\smc Problem.} Estimate  $\sharp_{inf}[c]$ in terms of the combinatorics
of $c$.
\vskip 2truemm

The problem itself is apparently motivated by the following classical
result due to M\"obius, see e.g. \cite {Ar3}.
\vskip 2truemm

{\smc Theorem.} Any embedded noncontractible curve on $\Bbb{RP}^2$ has at
least 3 inflection points.
\vskip 2truemm

The present paper contains some answers for the case when $c$ is a tree-like
curve, i.e. satisfies the condition that if $p$ is any double point of $c$
then $c\setminus p$ has  2   connected components. We plan to drop the
restrictive
assumption of tree-likeness in our next paper, see \cite {Sh}.
Classes of tree-like curves are naturally
enumerated by partially directed trees with a simple additional restriction
on directed edges, see \S 2. It was a pleasant surprise that for the classes of
tree-like curves there exists a (relatively) simple combinatorial criterion
characterizing  when $[c]$ contains a nonflattening curve, i.e.
$\sharp_{inf}[c]=0$ in terms of its tree. (Following V.Arnold we
use the word 'nonflattening' in this text as the synonym for the absence
of inflection points.) On the other hand, all attempts
to find a closed
formula for $\sharp_{inf}[c]$ in terms of partially directed trees failed.
Apparently such a formula does not exist, see the Concluding Remarks.
\vskip 2truemm

The paper is organized as follows. \S 2 contains some general information on
tree-like curves. \S 3 contains a criterion of noflattening. \S 4
presents some upper and lower bounds for $\sharp_{inf}[c]$. Finally, in
appendix we calculate the number of tree-like curves
having the same Gauss diagram.
\vskip 2truemm

For the general background on generic plane curves and their invariants the
author would recommend \cite {Ar1} and, especially, \cite {Ar2} which is an
excellent reading for a newcomer to the subject.
\vskip 2truemm

{\bf Acknowledgements.} It is a great pleasure to thank V.Arnold for the
formulation of the problem and his support and encouragement during many years
of our contact. His papers \cite {Ar1} and \cite {Ar2} in which he develops
Vassiliev-type theory for plane curves revived the whole area. It is difficult
 to overestimate the role of F.Aicardi for this project. Her suggestions have
essentially improved the original text and eliminated some mistakes in the
preliminary version.  Finally, the  author wants to acknowledge the hospitality
 and support of the Max-Planck-Institut during the preparation of this paper.
\bigskip

\heading \S 2. Some generalities on planar tree-like curves  \endheading

\bigskip
\centerline{\epsffile{fn1}}
\bigskip

Recall that a generic immersion $c:S^1\to \bR^2$ is called  {\it a
tree-like curve}
 if removing  any of its double points $p$ we get that $c\setminus p$ has
 2 connected components where $c$ also denotes the image set of $c$, see Fig.1.
Some of the  results below were first proved
in \cite {Ai} and later  found by the author independently. Recall that
{\it the Gauss diagram} of a generic immersion $c:S^1\to \bR^2$ is the original
circle $S^1$ with the set of all preimages of its double points, i.e. with the
set of all pairs $(\phi_1,\phi_2),...,(\phi_{2k-1},\phi_{2k})$ where
$\phi_{2j-1}$ and $\phi_{2j}$ are mapped to the same point on $\bR^2$ and
$k$ is the total number of double points of $c(S^1)$. One might think that
every pair of points $(\phi_{2j-1},\phi_{2j})$ is connected by an edge.

{\smc 2.1. Statement,} (see Proposition 2.1. in \cite {Ai}). A  generic
immersion $c:S^1\to \bR^2$
is a \tlc  iff its Gauss diagram is planar, i.e. can be drawn (including
edges connecting preimages of double points) on $\bR^2$
without selfintersections, see Fig.2.

\bigskip
\centerline{\epsffile{fn2}}
\bigskip

{\smc 2.2. Remark.} There is an obvious isomorphism between the set of all
planar Gauss
diagrams and the set of all planar connected trees.
 Namely, each planar Gauss diagram $GD$ corresponds to the following planar
tree.
Let us place a vertex in each connected component of $D^2\setminus GD$
 where $D^2$ is
the disc bounded
by the basic circle of  $GD$ and connect by edges
all vertices lying in the   neighboring connected components. The resulting
planar
tree is denoted by $Tr(GD)$. Leaves of $Tr(GD)$ correspond to the connected
components
with one neighbor.
On $GD$ these connected components  inherit  the natural cyclic order
according to their
position along the basic circle of $GD$. This cyclic order induces a natural
cyclic order on the set $Lv$ of leaves of its  planar tree $Tr(GD)$.
\vskip 2truemm

\subheading {Decomposition of a tree-like curve} Given a \tlc  $c: S^1\to
\bR^2$ we
decompose its image into a union of curvilinear polygons bounding
contractable
domains as follows. Take the planar Gauss diagram $GD(c)$ of $c$ and consider
the connected components of $D^2\setminus GD(c)$. Each such component
has a part of its boundary lying on $S^1$.
\vskip 2truemm

{\smc 2.3. Definition. } The image of the part of the boundary of a connected
component in
$D^2\setminus GD(c)$ lying on $S^1$ forms a closed nonselfintersecting
piecewise smooth
 curve (a curvilinear polygon) called {\it a building block} of $c$, see
Fig.3. (We
call vertices and edges
of building blocks {\it corners} and {\it sides} to distinguish them from
vertices and edges of planar trees used throughout the paper.)
\vskip 2truemm

The union of all building blocks constitutes the whole
\tl curve. Two building blocks have at most one common corner. If they have
a common corner then they are called {\it neighboring}.
\vskip 2truemm

{\smc 2.4. Definition.} A tree-like curve $c$ is called cooriented if every
its side is endowed with a coorientation, i.e. with a choice of local
connected component of the complement $\Bbb R^2\setminus c$ along the side.
(Coorientations of different sides are, in general, unrelated.) There are
two {\it continuous} coorientations of $c$ obtained by choosing one of two
possible coorientations of some side and extending it by continuity, see Fig.3.
\vskip 2truemm

\bigskip
\centerline{\epsffile{fn3}}
\bigskip

{\smc 2.5. Lemma.} Given a continuos coorientation of a \tlc  $\;c$ one
gets that all sides of any building block are either inward or outward
cooriented w.r.t. the interior of the block. (Since every building block
bounds a contractible domain  the outward and inward
coorientation have a clear meaning.)
\vskip 2truemm

{\smc Proof.} Simple induction on the number of building blocks.

\qed
\vskip 2truemm

{\smc 2.6. Definition. } Given a \tlc $c$ we associate to it the following
planar
partially directed tree $Tr(c)$. At first we take the undirected tree
$Tr(GD(c))$
where $GD(c)$ is the Gauss diagram of $c$, see Remark 2.2. (Vertices of
$Tr(GD(c))$ are in 1-1-correspondence with building blocks of $c$. Neighboring
blocks correspond to adjacent vertices of $Tr(GD(c))$.) For each pair
of neighboring building blocks $b_1$ and $b_2$ we do the following. If a
building
block $b_1$ contains a neighboring building block $b_2$ then we direct the
corresponding edge $(b_1,b_2)$ of $Tr(GD(c))$ from $b_1$ to $b_2$. The
resulting partially directed planar tree is denoted by $Tr(c)$.
Since $Tr(c)$ depends only on the class $[c]$ we will also use the notation
$Tr[c]$.
\vskip 2truemm

 We associate with each of two possible continuous   coorientations of $c$
 the following labeling of $Tr[c]$. For an  outward (resp. inward) cooriented
  building block we label by  $'+'$ (resp. $'-'$)
  the corresponding vertex of $Tr[c]$, see Fig.4.
\vskip 2truemm

{\smc 2.7. Definition.} Consider a partially directed tree $Tr$ (i.e. some of
its edges are directed).
$Tr$ is called {\it a noncolliding partially directed tree or ncpd-tree} if
no path of $Tr$
 contains edges pointing at each other.  The usual tree $Tr^\prime$
obtained by forgetting
directions of all edges
of $Tr$ is called {\it underlying}.
\vskip 2truemm

{\smc 2.8. Lemma.} a) For any \tlc $c$ its $Tr(c)$ is noncolliding;
\vskip 2truemm

b) the set of classes of nonoriented \tl curves is in 1-1-correspondence
with the set of all ncpd-trees on the nonoriented $\bR^2$.
\vskip 2truemm

{\smc Proof.} A connected component of \tl curves with a given Gauss
diagram is
uniquely determined by the enclosure of neighboring building blocks. The
obvious
restriction that if two building blocks contain the third one then one of
them is
contained in the other is equivalent to the noncolliding property.
(See an example in Fig. 4.)

\qed

\bigskip
\centerline{\epsffile{fn4}}
\bigskip

{\smc 2.9. Remark.} In terms of the above ncpd-tree one can easily describe
the Whitney index
(or the total rotation) of a given \tlc $c$ as well as the coorientation
of its building blocks.
 Namely, fixing the inward or outward coorientation of some building block
we determine the coorientation of any other building block as follows. Take
the (only)
path connecting the vertex corresponding to the fixed block with the vertex
corresponding to the other block. If the number of undirected edges in this
path is
odd then the coorientation changes and if this number is even then it is
preserved. (In other words, $Coor(b_1)=(-1)^{q(b_1,b_2)}Coor(b_2) $
where $q(b_1,b_2)$ is the number of undirected edges on the above path.)
\vskip 2truemm

{\smc 2.10. Lemma}, (see theorem 3.1 of \cite {Ai}).
 $$ind(c)=\sum_{b_i\in Tr(c)} Coor (b_i).$$
{\smc Proof.} Obvious.

\qed
\bigskip

\heading \S 3. Nonflattening of tree-like curves \endheading
\bigskip

In this section we give a criterion for nonflattening of a tree-like curve
in terms of its ncpd-tree, i.e. characterize all cases when
$\sharp_{inf}[c]=0$. (The author is aware of the fact that some of the
proofs below are rather sloppy since
they are based on very simple explicit geometric constructions on $\bR^2$ which
are not so easy to describe with complete rigor.)
\vskip 2truemm

{\smc 3.1. Definition.\;} A \tlc $c$ (or its class $[c]$) is called {\it
nonflattening} if
$[c]$
contains  a  generic immersion without inflection points.
\vskip 2truemm

{\smc 3.2. Definition.\;} {\it The  convex coorientation} of an arc $A$, image
of a smooth embedding $(0,1)\to \bR^2$ without inflection points is defined
as follows. The tangent line at any point $p\in A$ divides $\bR^2$ into two
parts. We choose at $p$ a vector transverse to $A$ and belonging to the
halfplane not containing $A$.
{\it The  convex coorientation} of a nonflattening  \tlc $c:S^1\to \bR^2$
is the convex coorientation of its arbitrary side extended
 by continuity to the whole curve, see Fig.5.

\bigskip
\centerline{\epsffile{fn5}}
\bigskip

{\smc 3.3. Definition.} Given a building block $b$ of a tree-like curve $c$
 we say that a corner $v$ of $b$
is  {\it of $\vee$-type (of $\wedge$-type resp.)} if  the interior angle
between the tangents to its sides at $v$ is bigger (smaller resp.) than
$180^0$, see Fig.6. (The interior
angle is the one contained in the interior of $b$.)
\vskip 2truemm

{\smc 3.4. Remark.\;} If $v$ is a $\vee$-type corner then the neighboring block
$b^\pr$ sharing
the corner $v$ with $b$ lies inside $b$, i.e. $\vee$-type corners are in
1-1-correspondence
with edges of the ncpd-tree of $c$ directed from the vertex corresponding
to $b$.
\vskip 2truemm

\bigskip
\centerline{\epsffile{fn6}}
\bigskip

{\smc 3.5. Criterion for nonflattening.} A \tlc $c$ is
nonflattening iff
the following 3 conditions hold for one of two possible coorientations of
its ncpd-tree, (see lemma 2.5).
\vskip 2truemm

a) all building $1$-gons (i.e. building blocks with one side) are
outward cooriented or, in terms of its tree,
all vertices of degree $1$ of $Tr[c]$ are labeled with $'+'$;
\vskip 2truemm

b) all building $2$-gons are outward cooriented or, in terms of its tree,
all vertices of degree $2$ are labeled with $'+'$;
\vskip 2truemm

c) the interior of any concave building block (=all sides are concave)
 with $k\ge 3$ sides contains at most $k-3$ other neigboring blocks or, in
terms of its tree,
any  vertex labeled by $'-'$ of degree $k\ge 3$ has at most $k-3$
leaving edges (i.e. edges directed from this vertex).
\medskip

{\smc Proof.} The necessity of a) - c) is rather obvious. Indeed, in the cases
  a) and b) a vertex of degree $\leq 2$ corresponds to a building block
with at most 2 corners. If such  a building block  belongs to a
nonflattening \tlc then it must be
globally convex and therefore outward cooriented w.r.t. the above convex
coorientation.
  For c) consider an inward cooriented (w.r.t. convex coorientation) building
block $b$ of a nonflattening curve. Such $b$ is a curvilinear polygon with
locally concave edges. Assuming
that $b$ has $k$ corners one gets that the sum of its interior angles is less
than $\pi(k-3)$. Therefore the number of $\vee$-type corners (or leaving
edges at the corresponding vertex) is less than $k-3$. (See Fig.7 for violations
of conditions a)-c). )
\vskip 2truemm

 Sufficiency of a)-c) is proved by a
relatively explicit construction. Given an ncpd-tree satisfying a) - c) let us
construct  a nonflattening curve with this tree using induction of the number
of vertices. (This will imply the sufficiency according to lemma 2.8.b.)
While constructing this curve inductively we provide additionally
that every building block is star-shaped with respect to some interior
point, i.e. the  segment connecting this point with a point on the boundary of
the block always lies in thew domain bounded by the block.
\vskip 2truemm

Case 1. Assume that the ncpd-tree contains an outward cooriented leaf
connected to an outward
cooriented vertex (and therefore the connecting edge is directed). Obviously,
the tree obtained by removal of this leaf is also a
ncpd-tree.  Using our induction we can construct a nonflattening curve
corresponding to the reduced tree and then, depending on the orientation of the
removed edge, either glue inside the appropriate locally convex building block a
small convex loop  (which is obviously possible) or glue a big locally
convex loop containing the whole curve. The possibility
to glue a big locally convex loop containing the whole curve is proved
independently in lemma 3.9 and corollary 3.10.
\vskip 2truemm

Case 2. Assume that all leaves are connected to inward cooriented vertices. (By
conditions a) and b)
these vertices are of degree $\ge 3$.) Using the ncpd-tree we can find
at least 1 inward cooriented vertex $b$ which is not smaller than any other
vertex, i.e.
the corresponding building block contains at least 1 exterior side. Let $k$ be
the degree of $b$ and $e_1,...,e_k$ be its edges in the cyclic order. (Each
$e_i$
is either undirected or leaving.)  By assumption c) the number of leaving edges
is at most $k-3$.
If we remove $b$ with all its edges then the remaining forest consists of
$k$ trees.
Each of the trees connected to $b$ by an undirected edge is an ncpd-tree.
We make every tree connected
to $b$ by a leaving edge into an ncpd-tree by gluing the undirected edge
instead of the
removed  directed and we mark the extra vertex we get.
By induction, we can construct $k$ nonflattening curves corresponding to
each of $k$ obtained
ncpd-trees. Finally, we have to glue
them to the corners of a locally concave $k$-gon with the sequence of $\vee$-
and $\wedge$-type
corners prescribed by $e_1,...,e_k$. The possibility of such a gluing
is proved independently in lemmas 3.11 and 3.12.

\qed

\bigskip
\centerline {\epsffile {fn7}}
\bigskip

 {\smc 3.6. An important construction.} The following operation called {\it
contracting homothety} will be extensively used below. It does not change the
class of a \tl curve and the number of inflection points.
\vskip 2truemm

Taking a  \tlc $c$ and some of its double points $p$ we
split $c\setminus p$ into 2 parts $c^+$ and $c^-$ intersecting only at $p$. Let
${\Omega}^+$ and
${\Omega}^-$ denote the domain bounded by the union of building blocks
 contained in  $c^+$ and $c^-$ resp. There are 2
options a) one of these unions
contains the other, say, ${\Omega}^-\subset {\Omega}^+$;
or b) $\overline{\Omega}^-\cap\overline{\Omega}^+=\{p\}$.
\vskip 2truemm

Let us choose some small neighborhood $\epsilon_p$ of the double point $p$
such that $c$ cuts $\epsilon_p$ in exactly four parts.
\vskip 2truemm

A contracting homothety with the centre $p$ in case a) consists of the
usual homothety $H$ applied to $c^-$ which places $H(c^-)$ into
$\epsilon_p$ followed
by smoothing   of the 2nd and higher derivatives of the union
$H(c^-)\bigcup c^+$
 at $p$. (This is always possible by changing $H(c^-)\bigcup c^+$ in some
even smaller neighborhood of $p$.) The resulting curve
 $c_1$ has the same ncpd-tree as $c$ and therefore belongs to $[c]$.
(Note that we do not construct an isotopy of $c$ and $c_1$ by applying
a family of homotheties with the scaling coefficient varying from $1$
to some small number. It suffices that the resulting curve $c_1$ has the
same ncpd-tree and therefore is isotopic to $c$.)
\vskip 2truemm

 In case b) we can apply a contracting homothety to either of 2 parts and
get 2 nonflattening \tl curves $c_1$ and $c_2$ isotopic to $c$  and such that
 either  $c^+=c_1^+$ while $c_1^-$ lies in an arbitrary small neighborhood
of $p$ or $c^-=c_2^-$ while
$c_2^+$ lies in an arbitrary small neighborhood of $p$. See Fig.8 for the
illustration of contracting homotheties.
\vskip 2truemm

{\smc 3.7. Definition.} Consider a bounded domain $\Omega$ in $\bR^2$ with
a locally strictly convex piecewise $C^2$-smooth boundary $\part \Omega$.
$\Omega$ is called {\it rosette-shaped} if for any side $e$ of $\part
\Omega$ there exists a point
$p(e)\in e$ such that $\Omega$ lies in one of the closed halfspaces
$\bR^2\setminus l_p(e)$
 w.r.t the tangent line $l_p(e)$ to $\part\Omega$ at $p(e)$.
\vskip 2truemm

{\smc 3.8. Remark.} For a  rosette-shaped $\Omega$ there exists a smooth convex
curve $\ga(e)$
containing $\Omega$ in its convex hull and tangent to $\part\Omega$ at exactly
one point
lying on a given side $e$ of $\part\Omega$.

\bigskip
\centerline {\epsffile {fn8}}
\bigskip

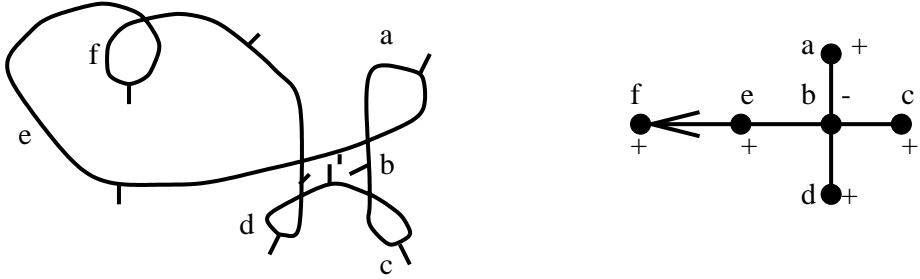
{\smc 3.9. Lemma.} Consider a nonflattening \tl curve $c$ with the convex
coorientation
and a locally convex building  block $b$ of $c$ containing at least one exterior
side, i.e. a side bounding the noncompact exterior domain on $\bR^2$. There
exists
a nonflattening curve $c^\prime$ isotopic to $c$ such that its building
block $b^\pr$
corresponding to $b$ bounds a rosette-shaped domain.
\vskip 2truemm

{\smc Proof. } Step 1. Let $k$ denote the number of corners of $b$.
Consider the connected components $c_1,...,c_k$ of $c\setminus b$.
By assumption that $b$ contains an exterior side one has that every $c_i$
lies either
inside or outside
 $b$ (can not contain $b$). Therefore using a suitable contracting homothety
 we can make every
$c_i$ small and lying in a prescribed small neighborhood of its corner
preserving
the nonflattening
property.
\vskip 2truemm

Step 2.  Take the standard  unit circle $S^1\subset
\bR^2$ and choose $k$ points on $S^1$. Then deform $S^1$ slightly into a
piecewise smooth locally
convex curve $\tilde{S^1}$ with the same sequence of $\vee$- and $\wedge$-type
corners as on $b$.
Now glue the small components $c_1,..., c_k$  (after an appropriate affine
transformation applied
to each $c_i$) to $\tilde{S^1}$ in the same order as they sit on $b$. The
resulting
curve $c^\pr$ is a nonflattening \tlc with the same ncpd-tree as $c$.

\qed
\vskip 2truemm

{\smc 3.10. Corollary.\;} Using the remark 3.8  one can glue a big
locally convex
loop containing the whole $c^\prime$ and tangent to $c^\pr$ at one point on
any exterior
edge and then deform this point of tangency into a double point and
therefore get the
nondegenerate \tl curve required in case 1.
\vskip 2truemm

Now we prove the supporting lemmas for case 2 of criterion 3.5.
\vskip 2truemm

Take any polygon $Pol$ with $k$ vertices and with the same sequence of $\vee$-
and $\wedge$-type
vertices as given by $e_1,...,e_k$; (see notations in the proof of case 2).
The existence of
such a polygon is exactly guaranteed by condition c), i.e. $k\ge 3$ and the
number of interior
angles $>\pi$ is less or equal than $k-3$. Deform it slightly to make it
into a locally concave
curvilinear polygon which we denote by $\widetilde{Pol}$.
\vskip 2truemm

{\smc 3.11. Lemma.} It is possible to glue a nonflattening curve $\tilde c$
 (after an appropriate diffeomorphism)
through its convex exterior edge to any $\wedge$-type vertex of
$\widetilde{Pol}$
placing it outside $\widetilde{Pol}$ and
preserving nonflattening of the union.
\vskip 2truemm

{\smc Proof. } We assume that the building block containing the side $e$ of
the curve $\tilde c$
to which we have to glue $v$ is rosette-shaped. We choose a point $p$ on
$e$ and substitute $e$ by 2 convex sides meeting transversally at $p$. Then we
apply to $\tilde c$ a linear transformation having its origin at $v$ in order
to a) make $\tilde c$ small, and b) to make the angle between
the new sides equal to the angle at the $\wedge$-type vertex $v$
to which we have to glue $\tilde c$. After that we glue $\tilde c$
by matching $v$ and $p$, making the tangent lines of $\tilde c$ at $p$
coinciding with the corresponding tangent lines of $\widetilde{Pol}$ at $p$
and smoothing the higher derivatives.

\qed
\vskip 2truemm

{\smc 3.12. Lemma.} It is possible to glue a nonflattening curve $\tilde c$
 after cutting off its exterior building
block with 1 corner to a $\vee$-type vertex of $\widetilde{Pol}$
and preserving the nonflattening of the union. The curve is placed inside
$\widetilde{Pol}$.
\vskip 2truemm

{\smc Proof. }
The argument is essentially the same as above. We cut off a convex exterior
loop  from $\tilde c$ and apply to the remaining curve a linear
transformation making it small and making the angle between 2 sides at the
corner where we have cut off a loop equal to the angle at the $\vee$-type
vertex.
Then we glue the result to the $\vee$-type vertex and smoothen the higher
derivatives.

\qed
\bigskip

\heading \S 4. {Upper and lower bounds of $\sharp_{inf}[c]$ for \tl curves}
\endheading
\bigskip

Violation of any of the above 3  conditions of
 nonflattening leads to the appearance of inflection points on a \tl curve
which are unremovable by diffeomorphisms.of $\Bbb R^2$. At first we reduce
the question about
the minimal number  $\sharp_{inf}[c]$ of inflection points of a class of
 \tl curves to a purely combinatorial
problem and then we shall give some upper and lower bounds for this number.
Some of the
geometric proofs are only sketched for the same reasons as in the previous
section. Since we are interested in inflections which survive under
the action of diffeomorphisms of $\bR^2$  we will assume from now on that
all considered curves have only locally unremovable inflection points, i.e.
those which do not disappear under arbitrarily small deformations of curves.
(For example, the germ $(t,t^4)$ is not interesting since its inflection
disappears after a arbitrarily small  deformation of the germ. One can
assume that the tangent line at every inflection point of any curve
we consider intersects the curve with the multiplicity 3.)
\vskip 2truemm

{\smc 4.1. Definition.} A generic immersion $c:S^1\to \bR^2$ the inflection
points of which coincide with some of its double points is called {\it
normalized.}
\vskip 2truemm

{\smc 4.2. Proposition.} Every \tl curve is isotopic to a normalized \tl
curve with at most the same number of inflection points.
\vskip 2truemm

{\smc Proof.}
Step 1. The idea of the proof is to separate building blocks as much as
possible and then substitute every block by a curvilinear polygon with
nonflattening sides.
Namely, given a \tl $c$ let us partially order the vertices of its ncpd-tree
$Tr[c]$ by choosing one
vertex as the root (vertex of level 1). Then we assign to all its adjacent
vertices level 2, etc.
The only requirement for the choice of the root is that all the directed edges
 point from the lower level to the higher. One can immediately see that
noncolliding
property garantees the existence of at least one root. Given such a partial
order
we apply consecutively a series of contracting homotheties to all double
points as
follows.
We start with double points which are the corners of the building block $b$
corresponding
to the root.
Then we apply contracting homotheties to all connected components of
$c\setminus b$ placing them into the  prescribed small neighborhoods of the
corners of $b$. Then we apply contracting homothety to all connected
components of $c\setminus (\cup b_i)$ where $b_i$ has level less or equal
$2$ etc. (See an example on Fig.9.)
(Again we do not need to construct an explicit isotopy of the initial and
final curves as soon as we know that they have the same ncpd-tree and,
therefore, are isotopic.)
Note that every building block except for the root has its father to which it is
attached through a $\wedge$-type corner since the root contains an exterior
edge.
The resulting curve $\tilde c$ has the same
type and number
of inflection points as $c$ and every building  lying in a prescribed
small neighborhood of the corresponding corner of its
father which does not intersect with other building blocks.

\bigskip
\centerline {\epsffile {fn9}}
\bigskip

 Step 2. Now we  substitute every side of every building block by a
nonflattening
arc not increasing the number of inflection points. Fixing some orientation
of $\tilde c$
we  assign at every double point 2 oriented tangent elements to 2 branches
of $\tilde c$
in the obvious way. Note that we can assume (after applying an
arbitrarily small deformation of $\tilde c$) that any 2 of these tangent
elements not sharing
the same vertex are in general position, i.e. the line connecting the
footpoints of the tangent
elements is different from both tangent lines.
\vskip 2truemm

Initial change. At first we will substitute every building block
of the highest level (which is necessarily a $1$-gon) by a convex loop.
According to our smallness assumptions there exists a smooth
(except for the corner) convex
loop gluing which instead of the building block will make the whole new curve
 $C^1$-smooth and isotopic to $\tilde c$ in the class of $C^1$-smooth curves.
This convex loop lies on the definite side w.r.t. both tangent lines at the
double point. Note that if the original removed building block lies wrongly
w.r.t. one (both resp.) tangent lines then it has at least 1 (2 resp.)
inflection points. After constructing a $C^1$-smooth curve we change it
slightly in a small neighborhood of the double point in order to
 provide for each branch a) if the branch of $\tilde c$ changes convexity
at the double point then we produce a smooth inflection at the double point;
 b) if the branch
does not change the convexity then we make it smooth. The above remark
garantees that the total number of inflection points does not increase.
\vskip 2truemm

Typical change. Assume that all blocks of level $> i$ already have
nonflattening sides. Take any block $b$ of level $i$. By the choice of the
root it has a unique $\wedge$-type corner with its father. The block $b$ has
a definite sequence of its $\vee$- and $\wedge$-type corners starting with
the attachment corner and going around $b$ clockwise. We want to
cut off all connected components of $c\setminus b$ which have level $>i$
then substitute $b$ by a curvilinear polygon with
nonflattening sides and then glue  back the blocks we cut off. Let us draw
the usual
polygon $Pol$ with the same sequence of $\vee$- and $\wedge$-type vertices
as for $b$. Now we will deform its sides into convex and concave arcs depending
on the sides of the initial $b$. The tangent elements to the ends of some side
of $b$ can be in one of  2 typical normal or 2 typical abnormal positions
(up to orientation-preserving affine transformations of $\bR^2$), see Fig.10.

\bigskip
\centerline {\epsffile {fn10}}
\bigskip

If the position of the tangent elements is normal then we deform the
corresponding side of $Pol$ to get a nonflattening arc with the same
position of the tangent elements as for the initial side of $b$. If
the position
is abnormal then we deform the side of $Pol$ to get a nonflattening arc
which
has the same position w.r.t. the tangent element at the beginning as the
original side of $b$.
Analogous considerations as before show that after gluing everything back
and smoothing the total number of inflections will not increase.

\qed
\vskip 2truemm

{\smc 4.3. Definition.} {\it A local coorientation} of a generic immersion
$c:S^1\to \bR^2$ is a free coorientation of each side of $c$ (i.e.
every arc between double points) which is, in general, discontinuous at
the double points. {\it The convex coorientation} of a normalized curve is
 its local coorientation which coincides inside each (nonflattening) side
with the convex coorientation of this side, see 3.2.

Given a tree-like curve $c$ with some local coorientation
we want to understand when there
exists a normalized tree-like curve $c'$ isotopic to $c$
whose convex coorientation coincides with a given local coorientation of $c$.
The following proposition is closely connected with the criterion of
nonflattening from \S 2 answers this question.
\vskip 2truemm

{\smc 4.4. Proposition (realizability criterion for a locally cooriented
tree-like
curve).}
 There exists a \tln curve with a given convex coorientation if and only if
the following 3 conditions  hold
\vskip 2truemm

a) every building $1$-gon is outward cooriented;
\vskip 2truemm

b) at least one side of every building $2$-gon is outward cooriented;
\vskip 2truemm

c) if some building block on $k$ vertices has only inward cooriented sides
then the domain it bounds contains at most $k-3$ neighboring building blocks.
\vskip 2truemm

(Note that only condition b) is somewhat different from that of criterion 3.5.)
\vskip 2truemm

{\smc Sketch of proof.} The necessity of a)-c) is obvious. These conditions
guarantee that all building blocks can be constructed.
It is easy to see that they are,
in fact, sufficient. Realizing each building block by some curvilinear
polygon with nonflattening sides we can glue them together in a global
normalized tree-like curve. Namely, we start from some building block which
contains an exterior edge. Then we glue all its neighbors to its corners.
(In order to be able to glue  them we make them small and adjust the gluing
angles by appropriate linear transformations.) Finally, we smoothen higher
derivatives at
all corners and then proceed in the same way for all new corners.

\qed


{\smc Combinatorial setup.}
The above proposition 4.4 allows us to reformulate the question about
the minimal number of inflection points $\sharp_{inf}[c]$ for \tl curves
combinatorially.
\vskip 2truemm

{\smc 4.5. Definition.\;} A local coorientation of a given tree-like curve $c$
 is called {\it admissible} if it satisfies the conditions a)-c) of
proposition 4.4.
\vskip 2truemm

{\smc 4.6. Definition.\;} Two sides of $c$ are called {\it neighboring}
if they share the same vertex and their tangent lines at this vertex coincide.
 We say that two neighboring sides in a locally cooriented curve $c$
{\it create an inflection point} if their coorientations are opposite.
 For a given local coorientation $Cc$ of a curve $c$ let $\sharp_{inf}(Cc)$
 denote the total number of created inflection points.
\vskip 2truemm

{\smc 4.7. Proposition} (combinatorial reformulation). For a given \tl
curve $c$ one has
$$\sharp_{inf}[c]=\min \sharp_{inf}(Cc)$$
where the minimum is taken over the set of all
admissible local coorientations $Cc$  of a tree-like curve $c$.
\vskip 2truemm

{\smc Proof.\;} This is the direct corollary of propositions 4.2 and 4.4.
Namely, for every
\tl curve $\tilde c$ isotopic to $c$ there exists a normalized curve
$\tilde c^\prime$
with at most the same number of inflection points. The number of inflection
points
of $\tilde c^\prime$ coincides with that of its convex coorientation
$Cc^\prime$.
On the other side,
for every admissible local coorientation $Cc$ there exists a normalized curve
$c^\prime$ whose convex coorientation coincides with $Cc$.

\qed
\vskip 2truemm

 \subheading {A lower bound}
\vskip 2truemm

 A natural lower bound for $\sharp_{inf}[c]$ can
be obtained
in terms of the ncpd-tree $Tr[c]$. Choose any of two continuous
coorientation of $c$
and the corresponding coorientation of $Tr[c]$, see \S 2. All 1-sided building
blocks of $c$
(corresponding to the leaves of $Tr[c]$) have the natural cyclic order.
(This order coincides with the natural cyclic order on all leaves of
$Tr[c]$  according to their
position on the plane.)
\vskip 2truemm

{\smc 4.8. Definition.} A neighboring pair of 1-sided building blocks (or of
leaves on $Tr[c]$)
is called {\it reversing} if the continuous coorientations of these blocks are
different. Let $\sharp_{rev}[c]$
denote the total number of reversing neighboring pairs of building blocks.
\vskip 2truemm

Note that $\sharp_{rev}[c]$ is even and independent on the  choice of the
 continuous coorientation of $c$. Moreover, $\sharp_{rev}[c]$ depends only
on the class
$[c]$ and therefore we can use the above notation instead of $\sharp_{rev}(c)$.
\vskip 2truemm

{\smc 4.9. Proposition.\;} $\sharp_{rev}[c]\le \sharp_{inf}[c].$
\vskip 2truemm

{\smc Proof.} Pick a point $p_i$ in each of the $1$-sided building blocks $b_i$
such that the side is locally convex near $p_i$ w.r.t. the interior of $b_i$.
(Such a choice is obviously possible since $b_i$ has just one side.) The proof
is accomplished by the following simple observation.
\vskip 2truemm

 Take an immersed segment $\ga: [0,1]\to \bR^2$ such that
$\ga(0)$ and $\ga(1)$ are not inflection points and the total number of
inflection points on $\ga$ is finite. At each nonflattening point $p$ of
$\ga$ we can choose the convex coorientation, see \S 3, i.e. since the
tangent line to $\ga$ at $p$ belongs locally to one connected component of
$\bR^2\setminus \ga$ we  can choose a transversal vector pointing at
that halfspace. Let us denote the convex coorientation at $p$ by $n(p)$.
\vskip 2truemm

{\smc 4.10. Lemma.\;} Assume that we have fixed a global continuous
coorientation $Coor$ of $\ga$. If $Coor(0)=n(0)$ and $Coor(1)=n(1)$ then
$\ga$ contains an even number of locally unremovable inflections.
If $Coor(0)=n(0)$ and $Coor(1)$ is opposite to $n(1)$ then $\ga$ contains
an odd number of locally unremovable  inflections.
\vskip 2truemm

{\smc Proof.}  Recall that we have assumed that all our inflection points
are unremovable by local deformations of the curve. Therefore passing through
such an inflection point the convex coorientation changes to the opposite.

\qed
\vskip 2truemm

\subheading {An upper bound}
\vskip 2truemm

{\smc 4.11. Definition.} If the number of $1$-sided building blocks is bigger
than $2$ then each pair of neighboring $1$-sided blocks of $c$
(leaves of $Tr[c]$ resp.)
is joined by the unique segment of $c$ not containing other $1$-sided
building blocks. We call this segment {\it a connecting path}.
If a connecting path joins a pair of neighboring $1$-sided blocks
 with the opposite coorientations (i.e. one block is inward cooriented
and the other is outward cooriented w.r.t. the continuous coorientation of
$c$)
then it is called {\it a reversing connecting path}.
\vskip 2truemm

{\smc 4.12. Definition.} Let us call   a
nonextendable sequence of $2$-sided building blocks not contained in each
other {\it a joint} of a \tl curve.
(On the level of its ncpd-tree one gets a sequence of degree 2 vertices
connected by undirected edges.)
\vskip 2truemm

Every joint consists of 2 smooth intersecting segments of $c$  called
{\it threads}  belonging to 2 different connecting paths.
\vskip 2truemm

{\smc 4.13. Definition.}
For every nonreversing connecting path $\rho$ in  $c$ we determine
 {\it its standard local coorientation} as follows. First
 we coorient its first side (which is the side of a $1$-gon) outward and then
extend this coorientation by continuity. (Since $\rho$ is nonreversing its
final
side will be outward cooriented as well.)
\vskip 2truemm

{\smc 4.14. Definition.} A joint is called {\it suspicious} if either
\vskip 2truemm

a) both its threads lie on nonreversing paths and both sides of some
2-sided block from this joint are inward cooriented w.r.t. the above standard
local coorientation of nonreversing paths; or
\vskip 2truemm

b) one thread lies on a nonreversing path and there exists a block
belonging to
this joint such that its side  lying on the nonreversing path is inward
cooriented  (w.r.t. the standard local coorientation of nonreversing paths); or
\vskip 2truemm

c) both threads lie on reversing paths.
\vskip 2truemm

Let $\sharp_{jt}$ denote the total number of suspicious joints.
\vskip 2truemm

{\smc 4.15. Definition.} A building block with $k$ sides is called
{\it suspicious} if it satisfies the following two conditions
\vskip 2truemm

a) it contains at least $k-3$ other blocks, i.e. at least
$k-3$ edges are leaving the corresponding vertex of the tree;
\vskip 2truemm

b) all sides lying on nonreversing paths are inward cooriented w.r.t.
the standard local coorientations of these nonreversing paths.
\vskip 2truemm

Let $\sharp_{bl}$ denote the total number of suspicious blocks.
\vskip 2truemm

{\smc 4.16. Proposition.} $$\sharp_{inf}[c]\le \sharp_{rev}[c]+2(\sharp_{jt}+
\sharp_{bl}).$$
\vskip 2truemm

{\smc Proof.} According to the statement 4.7 for any \tl curve $c$ one has
$\sharp_{inf}[c]\le \sharp_{inf}(Cc)$  where $Cc$ is some admissible local
 coorientation of $c$. Let us show that there exists
an admissible local coorientation with at most $\sharp_{rev}[c]+
2(\sharp_{jt}+
\sharp_{bl})$ inflection points, see Def 4.6. First we fix the standard local
 coorientations of all nonreversing paths. Then for each reversing path
we choose
any local coorientation with exactly one  inflection point (i.e. one
discontinuity
of local coorientation on the reversing path) to get the necessary
outward coorientations of all $1$-sided blocks.
Now the local coorientation of the whole fat ncpd-tree is fixed but
it is not admissible, in general.
In order to make it admissible we have to provide that conditions b) and c) of
Proposition 4.4 are satisfied for at most $\sharp_{jt}$ suspicious joints
and at most $\sharp_{bl}$ suspicious blocks. To make the local coorientation of
each such suspicious joint or block admissible we need to introduce at most
two additional inflection points for every suspicious joint or block.
 Proposition follows.

\qed
\bigskip

\heading \S 5. Concluding remarks. \endheading
\bigskip

In spite of the fact that there exists a reasonable criterion for nonflattening
in the class of \tl curves in terms of their ncpd-trees
the author is convinced that there is no closed formula for $\sharp_{inf}[c]$.
 Combinatorial reformulation of 4.6 reduces the calculation of
$\sharp_{inf}[c]$ to
a rather complicated discrete optimization problem which hardly is expected
to have an answer in a simple closed form. (One can even make speculations
about the computational
complexity of the above optimization problem.)
\vskip 2truemm

The lower and upper bounds presented in \S 4 can be improved by using much
more complicated characteristics of an ncpd-tree. On the other side, both
of them are
 sharp on some subclasses of ncpd-trees. Since in these cases a closed formula
 has not been obtained  the author  did not try to get the best possible
estimations here.
\vskip 2truemm

At the moment the author is trying to extend the results of this note
to the case of all generic curves in $\Bbb R^2$, see \cite {Sh}.
\bigskip

\heading \S 6. Appendix. Counting \tl curves with a given Gauss diagram.
\endheading
\bigskip

Combinatorial material of this section is not directly related to the
main content of the paper. It is a side product of the author's interest
in tree-like curves.
Here we calculate the number of different classes of
tree-like curves which have the same Gauss diagram.
\vskip 2truemm

{\smc 6.1. Proposition.} There exists a 1-1-correspondence between
classes of oriented tree-like curves with $n-1$ double points
on nonoriented $\bR^2$ and the set of all
planar ncpd-trees with $n$ vertices on oriented $\bR^2$.
\vskip 2truemm

{\smc Proof.} See \cite {Ai}.
\vskip 2truemm

{\smc 6.2. Definition.} For a given planar tree $Tr$ consider the subgroup
$Diff(Tr)$ of all
orientation-preserving diffeomorphisms of $\bR^2$ sending $Tr$
homeomorphically onto itself
as an embedded 1-complex.  The  subgroup $PAut(Tr)$ of the group $Aut(Tr)$ of
 automorphisms of $Tr$ as an abstract tree  induced by $Diff(Tr)$ is called
{\it the group of planar automorphisms of $Tr$}.
\vskip 2truemm

The following simple proposition gives a complete description of different
possible groups
$PAut(Tr)$. (Unfortunately, the author was unable to find a sutiable
reference for this but the proof is not too hard.)
\vskip 2truemm

{\smc 6.3. Statement.  \;}
\vskip 2truemm

\roster
\item  The group $PAut(Tr)$ of planar automorphisms of a given
planar tree $Tr$ is isomorphic to $\bZ/\bZ_p$ and is conjugate by an
appropriate diffeomorphism
to the rotation about some centre by multiples of $2\pi/p$.
\vskip 2truemm

\item If  $PAut(Tr)=\bZ/\bZ_p$ for $p>2$ then the above centre of rotation
is a vertex of $Tr$.
\vskip 2truemm

\item  For $p=2$  the centre of rotation is either a vertex of $Tr$ or the
middle of its edge.
\vskip 2truemm

\item  If the centre of rotation is a vertex of $Tr$ then the action
$PAut(Tr)$ on $Tr$
is free except for the centre and the quotient can be identified with a
connected
subtree $STr\subset Tr$ containing the centre.
\vskip 2truemm

\item  For $p=2$, if the centre is the middle of an edge, then the action of
$PAut(Tr)$ on $Tr$
 is free except for this edge.

\endroster
\vskip 2truemm

{\smc Sketch of proof.}
The action of $PAut(Tr)$ on the set $Lv(Tr)$ of leaves of $Tr$ preserves
the natural cyclic order
 on $Lv(Tr)$ and thus reduces to the $\bZ/\bZ_p$-action for some $p$. Now
each element $g\in
PAut(Tr)$ is determined by its action on $Lv(Tr)$ and thus the whole
$PAut(Tr)$ is isomorphic to
$\bZ/\bZ_p$. Indeed consider some $\bZ/\bZ_p$-orbit $O$ on $Lv(Tr)$ and all
vertices
of $Tr$ adjacent to
$O$. They are all pairwise different or all coincide since otherwise they
cannot form an orbit
of the action of diffeomorphisms on $Tr$.

\qed
\vskip 2truemm

{\smc 6.4. Proposition.\;} The number $\natural(GD)$ of all
classes of oriented \tl curves on nonoriented $\bR^2$ with a
given Gauss diagram $GD$ on $n$ vertices is equal
$$\natural(GD)=\cases a)\; 2^{n-1}+(n-1)2^{n-2}, \text { if } PAut(GD)
\text { is trivial; } \\
         b)\; 2^{2k-2}+(2k-1)2^{2k-3}+2^{k-2} \text { where  } n=2k,
         \text {if } PAut(GD)=\bZ/2\bZ \\
         \text {  and the  centre of rotation is the middle of the side; }\\
         c)\; 2^k +(2^{n-1}+(n-1)2^{n-2}-2^k)/p,\text { where } n=kp+1
         \text { and }\\  PAut(GD)=\bZ/p\bZ
         \text { for some prime p (including } PAut(GD)=\bZ/2\bZ\\
         \text { with a central vertex); }\\
         d)\; \text { for the general case see proposition 6.5 below. }
\endcases $$
\vskip 2truemm

{\smc Proof.\;} By Proposition 6.1 we enumerate ncpd-trees with a given
underlying planar tree $Tr(DG)$.
\vskip 2truemm

Case a).  Let us first calculate only ncpd-trees all edges of which are
directed. The
number of such ncpd-trees equals the number $n$ of vertices of $Tr(DG)$
since for any such tree
there exists such a source-vertex (all edges are directed from this vertex).
Now let us calculate
the number of ncpd-trees with $l$ undirected edges. Since $Aut(GD)$ is
trivial we can assume
that all vertices of $Tr(GD)$ are enumerated. There exist $\binom {n-1}{l}$
subgraphs in
$Tr(GD)$ containing $l$ edges and for each of these subgraphs there exist
$(n-l)$ ncpd-trees
with such a subgraph of undirected edges. Thus the total number
$\natural(GD)=\sum_{l=0}^{n-1}
\binom {n-1}{l}(n-l)=2^{n-1}+(n-1)2^{n-2}.$
\vskip 2truemm

Case b). The $\bZ/2\bZ$-action on the set of all ncpd-trees splits them
into 2 classes according
to the cardinality of orbits. The number
of $\bZ/2\bZ$-invariant ncpd-trees equals the number of all subtrees in a
tree on $k$ vertices
where $n=2k$ (since the source-vertex of such a tree necessarily lies in
the centre). The last
number equals $2^{k-1}$. This gives
$\natural(GD)=(2^{n-1}+(n-1)2^{n-2}-2^{k-1})/2+2^{k-1}=
2^{2k-2}+(2k-1)2^{2k-3}+2^{k-2}.$
\vskip 2truemm

Case c). The $\bZ/p\bZ$-action on the set of all ncpd-trees splits them
into 2 groups  according
to the cardinality of orbits. The number of $\bZ/p\bZ$-invariant ncpd-trees
equals the number
of all subtrees in a tree with $k+1$ vertices where $n=pk+1$ (since the
source-vertex of such a tree
 lies in the centre). The last number equals $2^{k}$. This gives $\natural(GD)=
2^{k}+(2^{n-1}+(n-1)2^{n-2}-2^{k})/p.$
\vskip 2truemm

{\smc 6.5. Proposition.} Consider a Gauss diagram $GD$ with a tree $Tr(GD)$
 having $n$ vertices
 which has $Aut(Tr)=\bZ/p\bZ$ where $p$ is not a prime. Then for each
nontrivial factor $d$ of
$p$ the number of ncpd-trees with $\bZ/d\bZ$ as their group of symmetry equals
$$\sum_{d^\prime\vert d}\mu(d^\prime) 2^\frac {kd}{d^\prime} $$
where $\mu(d^\prime)$ is the M\"obius function. (This gives a rather
unpleasant expression
for the number of all \tl curves with a given $GD$ if $p$ is an arbitrary
positive integer.)
\vskip 2truemm

{\smc Proof.\;} Consider for each $d$ such that $d\vert p$ a subtree
$STr_d$ with $km+1$
vertices $m=\frac {n-1}{d}$ 'spanning' $Tr$ with respect to the
$\bZ/d\bZ$-action. The number of ncpd-trees invariant
at least w.r.t $\bZ/d\bZ$ equals $2^{km}$ where $p=dm$ and $n=kp+1$. Thus
by the inclusion-exclusion
formula one gets that the number of ncpd-trees invariant exactly w.r.t.
$\bZ/d\bZ$ equals
$\sum_{d^\prime\vert d}\mu(d^\prime) 2^\frac {kd}{d^\prime}. $

\qed
\vskip 2truemm

{\smc Problem.} Calculate the number of ncpd-trees with a given underlying
tree and of a given index.

\enddocument